\title{ ~~\\ Asymptotically exact heuristics for prime divisors of
the sequence $\{a^k+b^k\}_{k=1}^{\infty}$}
\author{Pieter Moree}
\documentclass[12pt]{article}
\usepackage{amssymb, latexsym, amsfonts}
\def\@ptsize{2}
\newtheorem{Thm}{Theorem}
\newtheorem{Lem}{Lemma}
\newtheorem{Cor}{Corollary}

\newtheorem{Prop}{Proposition}
\newcommand{\qed}{\hfill $\Box$}

\begin{document}
\date{}
\maketitle 
{\def\thefootnote{}
{\def\thefootnote{}
}
\begin{abstract}
\noindent Let $N_{a,b}(x)$ count the number of primes $p\le x$ with $p$ 
dividing $a^k+b^k$ for some $k\ge 1$. It is known that 
$N_{a,b}(x)\sim c(a,b)x/\log x$ for some rational number $c(a,b)$ 
that depends in a rather intricate way on $a$ and $b$. A simple 
heuristic formula for $N_{a,b}(x)$ 
is proposed and it is proved that it is asymptotically exact, i.e.
has the same asymptotic behaviour as $N_{a,b}(x)$. Connections with
Ramanujan sums and character sums are discussed.
\end{abstract}
\section{Introduction}
Let $p$ be a prime (indeed, throughout this note the letter $p$ will be used to
indicate primes). Let $g$ be a non-zero rational number. By $\nu_p(g)$ we denote the exponent
of $p$ in the canonical factorisation of $g$. If $\nu_p(g)=0$, 
then by ord$_g(p)$ we denote the
smallest positive integer $k$ such that $g^k\equiv 1({\rm mod~}p)$. If $k=p-1$, then
$g$ is said to be a {\it primitive root} mod $p$. If $g$ is a primitive root mod $p$, 
then $g^j$ is a primitive root mod $p$ iff gcd$(j,p-1)=1$. There are thus
$\varphi(p-1)$
 primitive roots mod $p$ in $(\mathbb Z/p\mathbb Z)^*$, where $\varphi$ denotes Euler's
 totient function.\\
\indent Let  $\pi(x)$ denote the number of primes $p\le x$ and 
$\pi_g(x)$ the number of primes $p\le x$ such that $g$ is a primitive root mod $p$. 
Artin's celebrated primitive root conjecture (1927) states that if $g$ is an
integer with $|g|>1$ and $g$ is not a square, then for some positive rational number $c_g$
we have $\pi_g(x)\sim c_gA\pi(x)$, as $x$ tends to infinity. Here $A$ denotes
{\it Artin's constant} 
$$A=\prod_p\left(1-{1\over p(p-1)}\right)=0.3739558136\dots$$ 
Hooley \cite{Hooley}, under
assumption of the Generalized Riemann Hypothesis (GRH), established Artin's conjecture and
explicitly evaluated $c_g$.\\
\indent It is an old heuristic idea that the behaviour of $\pi_g(x)$ should be mimicked
by $H_1(x)=\sum_{p\le x}\varphi(p-1)/(p-1)$, the idea being that the `probability' that $g$ is
a primitive root mod $p$ equals $\varphi(p-1)/(p-1)$ (since this is the density of
primitive roots in $(\mathbb Z/p\mathbb Z)^*$). Using the Siegel-Walfisz theorem 
(see Lemma \ref{walvis} below), it is
not difficult to show, unconditionally, that  $H_1(x)\sim A\pi(x)$. Although
true for many $g$ and also on average, it is however not always true, under GRH,
that $\pi_g(x)\sim H_1(x)$, i.e., the heuristic $H_1(x)$ is
not always asymptotically exact. Nevertheless, 
Moree \cite{M1} found a quadratic modification, $H_2(x)$, of the above heuristic $H_1(x)$
involving the Legendre symbol that is always asymptotically exact (assuming GRH).\\
\indent A prime $p$ is said to divide a sequence $S$ of integers, if it divides at least
one term of the sequence $S$ (see \cite{Ballot} for
a nice introduction to this topic). Several authors studied the problem of characterising
(prime)divisors of the sequence $\{a^k+b^k\}_{k=1}^{\infty}$. Hasse \cite{Hasse} seems to have
been the first to consider the Dirichlet density of prime divisors of such
sequences. Later authors, e.g., Odoni \cite{Odoni}
and Wiertelak \cite{Wiertelak} strengthened the analytic aspects of his work.
The best result to date, in the formulation of \cite{M0}, seems to be as follows 
(recall that Li$(x)=\int_2^x dt/\log t$ is the logarithmic integral):
\begin{Thm}
\label{oud}
Let $a$ and $b$ be non-zero integers. Put $r=a/b$. Assume that $r\ne \pm 1$.
Let $\lambda$ be the largest integer such that $|r|=u^{2^{\lambda}}$, with
$u$ a rational number. Let $\varepsilon={\rm sign}(r)$ and
$L=\mathbb Q(\sqrt{u})$. We have
$$N_{a,b}(x)=\delta(r){\rm Li}(x)+O\left({x(\log \log x)^4\over \log^3 x}\right),$$
where the implied constant may depend on $a$ and $b$ and $\delta(r)$, a
rational number, is given
in Table {\rm 1}.
\end{Thm}
\medskip
\centerline{{\bf Table 1:} The value of $\delta(r)$}
\medskip
\medskip
\begin{center}
\begin{tabular}{|c|c|c|c|}
\hline
$L$  & $\lambda$ & $\varepsilon=+1$ & 
$\varepsilon=-1$ \\
\hline
$L\ne \mathbb Q(\sqrt{2})$  & $\lambda\ge 0$ & $2^{1-\lambda}/3$ & 
$1-2^{-\lambda}/3$ \\
\hline
$L=\mathbb Q(\sqrt{2})$  & $\lambda=0$ & $17/24$ & 
$17/24$ \\
\hline
$L=\mathbb Q(\sqrt{2})$  & $\lambda=1$ & $5/12$ & 
$2/3$ \\
\hline
$L=\mathbb Q(\sqrt{2})$  & $\lambda\ge 2$ & $2^{-\lambda}/3$ & 
$1-2^{-1-\lambda}/3$ \\
\hline
\end{tabular}
\end{center}
\medskip
\noindent Starting point in the proof of Theorem \ref{oud} is the observation that 
$p\nmid 2ab$ divides the sequence $\{a^k+b^k\}_{k=1}^{\infty}$ iff ord$_r(p)$ is 
even, where $r=a/b$. The condition that ord$_r(p)$ be even is weaker than the
condition that ord$_r(p)=p-1$ and now the analytic tools are strong enough to establish
an unconditional result.\\
\indent Note that $\delta(r)$ does not depend on $\epsilon$ in case $\lambda=0$. For a `generic'
choice  of $a$ and $b$, $L$ will be different from $\mathbb Q(\sqrt{2})$ and $\lambda$ will
be zero and hence $\delta(a/b)=2/3$. It is not difficult to show \cite{M4} that the average
density of elements of even order in a finite field of prime cardinality also equals $2/3$.\\
\indent In this note analogs of $H^{(1)}_{a,b}(x)$ and $H^{(2)}_{a,b}(x)$ of $H_1(x)$
and $H_2(x)$ will be introduced and it will be shown that $H^{(2)}_{a,b}(x)$ is
always asymptotically exact. This leads to the following main result 
(where $\pi(x;k,l)$ denotes the number of primes $p\le x$ satisfying $p\equiv l({\rm mod~}k)$ and $(*/p)$ 
denotes the Legendre symbol):
\begin{Thm}
\label{main}
Let $a$ and $b$ be non-negative natural numbers.
Let $N_{a,b}(x)$ count the number of primes $p\le x$ that divide some
term $a^k+b^k$ in the sequence $\{a^k+b^k\}_{k=1}^{\infty}$. Put
$r=a/b$ and $\epsilon={\rm sgn}(a/b)$. Assume that $r\ne \pm 1$. Let $h$ be the largest integer such
that $|r|=r_0^h$ for some $r_0\in \mathbb Q$ and $h\ge 1$.
Put $e=\nu_2(h)$. If $\epsilon=1$, then
$$N_{a,b}(x)=\pi(x;2^{e+1},1)-2^{e+1}\sum_{p\le x,~(r_0/p)=1\atop \nu_2(p-1)>e}2^{-\nu_2(p-1)}
+O\left({x(\log \log x)^4\over \log^3 x}\right),$$
and if $\epsilon=-1$, then
$$N_{a,b}(x)=\pi(x)-\sum_{p\le x,~(r_0/p)=-1\atop \nu_2(p-1)=e+1}1-2^{e+1}\sum_{p\le x,~(r_0/p)=1\atop 
\nu_2(p-1)>e+1}2^{-\nu_2(p-1)}
+O\left({x(\log \log x)^4\over \log^3 x}\right),$$
where the
implied constants depend at most on $a$ and $b$.
\end{Thm}

\section{Preliminaries} 
The proof of Theorem \ref{main} requires a result from analytic number theory: the 
Siegel-Walfisz theorem, see e.g., \cite[Satz 4.8.3]{Prachar}. For notational convenience
we write $(a,b)$ instead of gcd$(a,b)$.
\begin{Lem} 
\label{walvis} 
Let $C>0$ be arbitary. There exists $c_1>0$ such that
$$\pi(x;k,l)={{\rm Li}(x)\over \varphi(k)}+O(x e^{-c_1\sqrt{\log x}}),$$
uniformly for $1\le k\le \log^C x$, $(l,k)=1$, where the 
implied constant depends at most on $C$.
\end{Lem} 
Our two heuristics will be based on the following
elementary observation in group theory.
\begin{Lem}
\label{Eerste} $~$\\
\noindent {\rm 1)} Let $h\ge 1$ and $w\ge 0$ be integers. Let $G$ be a cyclic group of order $n$.
Let $G^h=\{g^h:g\in G\}$ and
$G_{w}^h=\{g^h:\nu_2({\rm ord}(g^h))=w\}$. We have
$\# G^h=n/(n,h)$ and $\# G_{0}^h=2^{-\nu_2(n/(n,h))}n/(n,h)$.
Furthermore, for $w\ge 1$, we have
\begin{equation}
\label{subtieler}
\#G_w^h=\cases{2^{w-1-\nu_2(n/(n,h))}n/(n,h) &if $\nu_2(n/(n,h))\ge w$;\cr
0 & otherwise.}
\end{equation}
{\rm 2)} If $\nu_2(h)\ge \nu_2(n)$, then every element in $G^h$ has odd order.
If $\nu_2(h)<\nu_2(n)$, then $G_0^h\subseteq G^{2h}$.\\
{\rm 3)} We have
$$G_1^{h}\subseteq \cases{G^h\backslash G^{2h} &if $\nu_2(n)=\nu_2(h)+1$;\cr
G^{2h} &if $\nu_2(n)>\nu_2(h)+1$.}$$
If $\nu_2(n)\le \nu_2(h)$, then $G_1^{h}$ is empty.
\end{Lem}
{\it Proof}. 1) Let $g_0$ be a generator of $G$. On noting that
$g_0^{m_1}=g_0^{m_2}$ iff $m_1\equiv m_2({\rm mod~}n)$, the 
proof becomes a 
simple exercise in solving linear congruences. In this way one infers
that $G^h=\{g_0^{hk}~:~1\le k\le n/(n,h)\}$ and hence
$\# G^h=n/(n,h)$. Note that ord$(g_0^{hk})$ is the smallest positive integer
$m$ such that $n/(n,h)$ divides $mk$. Thus ord$(g_0^{hk})$ will be odd iff
$\nu_2(k)\ge \nu_2(n/(n,h))$. 
Using this observation we obtain that
\begin{equation}
\label{gnul}
G_{0}^h=\{g_0^{hk}~:~1\le k\le {n\over (n,h)},~\nu_2(k)\ge 
\nu_2({n\over (n,h)})\}
\end{equation}
and hence $\#G_{0}^h=2^{-\nu_2(n/(n,h))}n/(n,h)$. Similarly
$$G_{w}^h=\{g_0^{hk}~:~1\le k\le {n\over (n,h)},~\nu_2(k)=\nu_2({n\over (n,h)})-w\}$$
and hence we obtain (\ref{subtieler}).\\
2) If $\nu_2(h)<\nu_2(h)$, then using (\ref{gnul}) we infer that
$$G_0^h\subseteq \{g_0^{hm}:1\le m\le {n\over (n,h)},~\nu_2(m)\ge 1\}
=\{g_0^{2hk}:1\le k\le {n\over (n,2h)}\}=G^{2h},$$
where we have written $m=2k$ and used that $(n,2h)=2(n,h)$.\\
3) Similar to that of part 2. \qed\\

\noindent Remark. Note that $G^h$ and $G_0^h$ with the induced group operation from
$G$ are actually subgroups of $G$.

\section{Two heuristic formulae for $N_{a,b}(x)$}
In this section we propose two heuristics for $N_{a,b}(x)$; one more
refined than the other. Starting point is the observation that
if $p\nmid 2ab$ divides the sequence $\{a^k+b^k\}_{k=1}^{\infty}$ if
and only if ord$_r(p)$ is even, where $r=a/b$. Let $h$
be the largest integer such that we can write $|r|=r_0^h$ with $r_0$ 
a rational number. Let $\epsilon={\rm sgn}(r)$.\\
\indent We will use Lemma \ref{Eerste} in 
the case $G=G_p:=(\mathbb Z/p\mathbb Z)^*\cong \mathbb F_p^*$. 
The first heuristic approximation we consider is
$$K_{a,b}^{(1)}(x)=\sum_{p\le x,~p\nmid 2ab}{\# G_{p,(1-\epsilon)/2}^h\over \# G_p^h},$$
where $K_{a,b}^{(1)}(x)$ is supposed to be an heuristic for the number of primes $p\le x$
such that ord$_r(p)$ is odd. From our results below it will follow that $\lim_{x\rightarrow \infty}
K_{a,b}^{(1)}(x)/\pi(x)$ exists. Note that in case $h=1$, this 
limit is the average density of elements
of odd order (if $\epsilon=1$), respectively of order congruent to $2({\rm mod~}4)$ (if
$\epsilon=-1$). For a more detailed investigation of the average number of elements having
order $\equiv a({\rm mod~}d)$ vide \cite{M4}.\\
\indent Suppose that $p\nmid 2ab$. By assumption $r\in \epsilon G_p^h$. 
In the case $\epsilon=1$, the latter set  has
$\# G_{p,0}^h$ elements having odd order and so, in some sense,
${\# G_{p,0}^h/\# G_p^h}$ is the probability that ord$_r(p)$ is
odd. This motivates the definition of $K_{a,b}^{(1)}(x)$ in case $\epsilon=1$.
In case $\epsilon=-1$ we use the observation that for $p$ is odd, $-r_0^h$ has
odd order iff $r_0^h$ has order congruent to $2({\rm mod~}4)$. Thus the
elements in $-G_p^h$ of odd order are precisely the elements having order 
$2({\rm mod~}4)$ in $G_p^h$ and hence have cardinality $\# G_{p,1}^h$. 
On using part 1 of Lemma \ref{Eerste} we infer that
$K_{a,b}^{(1)}(x)=\sum_{p\le x,~p\nmid 2ab}k_{a,b}^{(1)}(p)$ with
\begin{equation}
\label{lokaaleen}
k_{a,b}^{(1)}(p)=\cases{(1+\epsilon)/2 & if $\nu_2(p-1)\le e$;\cr
2^{e-\nu_2(p-1)} & if $\nu_2(p-1)>e$.}
\end{equation}
An
heuristic $H_{a,b}^{(1)}(x)$ for $N_{a,b}(x)$ is now obtained on merely setting
$H_{a,b}^{(1)}(x)=\pi(x)-K_{a,b}^{(1)}(x)$. 
Put $\omega(n)=\sum_{p|n}1$.
On using (\ref{lokaaleen}) we then
infer that
$$H_{a,b}^{(1)}(x)=\pi(x;2^{e+1},1)-2^{e}\sum_{p\le x\atop \nu_2(p-1)>e}2^{-\nu_2(p-1)}+
O(\omega(ab)).$$ 
if $\epsilon=1$
and $$H_{a,b}^{(1)}(x)=\pi(x)-2^{e}\sum_{p\le x\atop \nu_2(p-1)>e}2^{-\nu_2(p-1)}+O(\omega(ab)).$$ 
if $\epsilon=-1$.\\
\indent In the context of (near) primitive roots it
is known that the analoga of $H_{a,b}^{(1)}(x)$ do not always, assuming
GRH, exhibit the correct 
asymptotic behaviour, but that an appropriate `quadratic' heuristic, i.e.
an heuristic taking into account Legendre symbols, always has the
correct asymptotic behaviour \cite{M1, M2, M3} (in \cite{M3} 
the main result of \cite{M2} is proved in a different and much shorter way). With this in 
mind, we propose a second, more refined, heuristic: $H_{a,b}^{(2)}(x)$.\\
\indent If $\nu_p(r)=0$ we
can consider $|r|=r_0^h$ and $r_0$ as elements of $G_p$. 
We write $(r_0/p)=1$ if $r_0$ is a square in $G_p$ and $(r_0/p)=-1$ otherwise.\\
\indent First consider the case where $\epsilon:={\rm sgn}(r)=1$. If
$\nu_2(p-1)\le e:=\nu_2(h)$, then $r$ has odd order by part 2 of Lemma \ref{Eerste}.
If $\nu_2(p-1) > \nu_2(h)$ and $(r_0/p)=-1$, then $r\in G_p^h$, but
$r\not\in G_p^{2h}$ (by part 2 of Lemma \ref{Eerste} again). It then follows that $r$ has even order.
On the other hand, if $(r_0/p)=1$ then $r\in G_p^{2h}$. This suggests to take
$$K_{a,b}^{(2)}(x)=\sum_{p\le x,~\nu_2(p-1)\le e}1+\sum_{p\le x,~(r_0/p)=1\atop \nu_2(p-1)>e}{\# G_{p,0}^h
\over \# G_p^{2h}},$$ where furthermore we require that $p\nmid 2ab$.
A similar argument, now using part 3 instead of part 2 of Lemma \ref{Eerste}, leads to the choice
$$K_{a,b}^{(2)}(x)=\sum_{p\le x,~(r_0/p)=-1\atop \nu_2(p-1)=e+1}{\# G_{p,1}^h
\over \# G_p^{2h}}+\sum_{p\le x,~(r_0/p)=1\atop \nu_2(p-1)>e+1}{\# G_{p,1}^h
\over \# G_p^{2h}},$$
in case $\epsilon=-1$, where again we furthermore require 
that $p\nmid 2ab$. We obtain $K_{a,b}^{(2)}(x)=\sum_{p\le x,~p\nmid 2ab}k_{a,b}^{(2)}(p)$, with
\begin{equation}
\label{lokaaltwee}
k_{a,b}^{(2)}(p)=\cases{(1+\epsilon)/2 & if $\nu_2(p-1)\le e$;\cr
(1+\epsilon ({r_0\over p}))/2 & if $\nu_2(p-1)=e+1$;\cr
(1+({r_0\over p}))2^{e-\nu_2(p-1)} & if $\nu_2(p-1)>e+1$.}
\end{equation}
Now we put $H_{a,b}^{(2)}(x)=\pi(x)-K_{a,b}^{(2)}(x)$ as before.
On invoking Lemma \ref{Eerste}, $H_{a,b}^{(2)}(x)$ can then be more explicitly written as
\begin{equation}
\label{heurieen} 
H_{a,b}^{(2)}(x)=\pi(x;2^{e+1},1)-2^{e+1}\sum_{p\le x,~(r_0/p)=1\atop \nu_2(p-1)>e}2^{-\nu_2(p-1)}
+O(\omega(ab)),
\end{equation}
if $\epsilon=1$ and
\begin{equation}
\label{heuritwee}
H_{a,b}^{(2)}(x)=\pi(x)-
\sum_{p\le x,~(r_0/p)=-1\atop \nu_2(p-1)=e+1}1-
2^{e+1}\sum_{p\le x,~(r_0/p)=1\atop \nu_2(p-1)>e+1}2^{-\nu_2(p-1)}+O(\omega(ab)),
\end{equation}
if $\epsilon=-1$.

\section{Asymptotic analysis of the heuristic formulae}
In this section we determine the asymptotic behaviour of $H_{a,b}^{(1)}(x)$ and
$H_{a,b}^{(2)}(x)$. 
We adopt the notation from Theorem \ref{main} and in addition write $D$ for the
discriminant of $\mathbb Q(\sqrt{r_0})$. Note that $D>0$.
\begin{Thm} 
\label{maintwee}
Let $A>0$ be arbitrary. The implied constants below depend at most on $A$.\\
{\rm 1)} We have
$$H_{a,b}^{(1)}(x)=\delta_1(r){\rm Li}(x)+O(x\log^{-A}x)+O(\omega(ab)),$$
where $$\delta_1(r)=\cases{2^{1-e}/3 & if $\epsilon=+1$;\cr
1-2^{-e}/3 & if $\epsilon=-1$.}$$
In particular, if $L\ne \mathbb Q(\sqrt{2})$, then $H_{a,b}^{(1)}(x)$ is
an asymptotically exact heuristic for $N_{a,b}(x)$.\\
{\rm 2)}  We have
$$H_{a,b}^{(2)}(x)=\delta(r){\rm Li}(x)+O(D^2x\log^{-A}x)+O(\omega(ab)).$$
In particular, $H_{a,b}^{(2)}(x)$ is
an asymptotically exact heuristic for $N_{a,b}(x)$.
\end{Thm}
The proof of part 2 requires a few facts from algebraic number theory, the proof of part 1 does
not even require that and is an easier variant of the proof of part 2 (and is
left to the interested reader). The proof of part 2 rests on a few lemmas.
\begin{Lem}
\label{voorgaande}
Let $n$ be a non-zero integer and $K=\mathbb Q(\sqrt{n})$ a quadratic number field of discriminant $\Delta$.
Let $A>1$ and $C>0$ be positive real numbers. Then
$$\sum_{p\le x,~(n/p)=1\atop \nu_2(p-1)=k}1={\rm Li}(x)\left({1\over [K(\zeta_{2^k}):\mathbb Q]}
-{1\over [K(\zeta_{2^{k+1}}):\mathbb Q]}\right)+O\left({|\Delta|x\over \log^A x}\right),$$
uniformly in $k$ with $k$ satisfying $2^{k+3}|\Delta|\le \log^C x$, where the implied constant depends
at most on $A$ and $C$.
\end{Lem}
{\it Proof}. By quadratic reciprocity a prime $p$ satisfies $(n/p)=1$ iff $p$ is in a certain set of 
congruences classes modulo $4|\Delta|$. Thus the primes we are counting in our sum are precisely
the primes that belong to certain congruences classes modulo $2^{k+2}|\Delta|$, but do not belong
to certain congruence classes of modulus $2^{k+3}|\Delta|$. The total number of congruence classes
involved is less than $8|\Delta|$. Now apply Lemma \ref{walvis}. This yields the result but with
an, as yet, unknown density.\\
\indent On the other hand, the primes $p$ that are counted are precisely the primes $p\le x$ that split
completely in the normal number field $K(\zeta_{2^k})$, but do not split
completely in the normal number field
$K(\zeta_{2^{k+1}})$. If $M$ is any normal extension then it is a consequence of Chebotarev's density
theorem that the set of primes that split completely in $M$ has density $1/[M:\mathbb Q]$. On using
this, the proof is completed. \qed 

\begin{Lem}
Let $m$ be fixed. With the notation as in the previous lemma we have
$$\sum_{p\le x,~(n/p)=1\atop \nu_2(p-1)\ge m}2^{-\nu_2(p-1)}={\rm Li}(x)\sum_{k=m}^{\infty}
{1\over 2^k}\left({1\over [K(\zeta_{2^k}):\mathbb Q]}
-{1\over [K(\zeta_{2^{k+1}}):\mathbb Q]}\right)+O\left({\Delta^2x\over \log^A x}\right),$$
where the implied constant depends at most on $A$.
\end{Lem}
{\it Proof}. We have
$$\sum_{p\le x,~(n/p)=1\atop \nu_2(p-1)\ge m}2^{-\nu_2(p-1)}=\sum_{k=m}^{m_1}
\sum_{p\le x,~(n/p)=1\atop \nu_2(p-1)=k}2^{-k}+O({x\over 4^{m_1}}),$$
where we used the trivial bound $\sum_{p\le x,~\nu_2(p-1)\ge m_1}2^{-\nu_2(p-1)}=O({x/4^{m_1}})$.
Choose $m_1$ to be the largest integer such that $2^{m_1+3}|\Delta|\le \log^C x$. Apply
Lemma \ref{voorgaande} with any $C>A/2$. It follows that
\begin{eqnarray}
\sum_{p\le x,~(n/p)=1\atop \nu_2(p-1)\ge m}2^{-\nu_2(p-1)}&=& {\rm Li}(x)\sum_{k=m}^{m_1}
{1\over 2^k}\left({1\over [K(\zeta_{2^k}):\mathbb Q]}
-{1\over [K(\zeta_{2^{k+1}}):\mathbb Q]}\right)+O({x\over 4^{m_1}});\nonumber\cr
&=& {\rm Li}(x)\sum_{k=m}^{\infty}
{1\over 2^k}\left({1\over [K(\zeta_{2^k}):\mathbb Q]}-
{1\over [K(\zeta_{2^{k+1}}):\mathbb Q]}\right)+O({x\over 4^{m_1}}),\nonumber
\end{eqnarray}
where we used that $\varphi(2^k)\le [K(\zeta_{2^k}):\mathbb Q]\le 2\varphi(2^k)$. 
On noting that $O(x/4^{m_1})=O(\Delta^2x\log^{-A}x)$, the result follows.\qed

\begin{Lem}
\label{bijna}
We have $H_{a,b}^{(2)}(x)=\delta_2(r){\rm Li}(x)+O(D^2x\log^{-A}x)+O(\omega(ab))$, where
\begin{equation}
\label{pluseen}
\delta_2(r)={1\over 2^e}-2^{e+1}\sum_{k=e+1}^{\infty}{1\over 2^k}\left({1\over [L(\zeta_{2^k}):\mathbb Q]}
-{1\over [L(\zeta_{2^{k+1}}):\mathbb Q]}\right)
\end{equation}
if $\epsilon=1$ and
$$\delta_2(r)=1-{1\over 2^{e+1}}
+{1\over [L(\zeta_{2^{e+1}}):\mathbb Q]}-{1\over [L(\zeta_{2^{e+2}}):\mathbb Q]}$$
\begin{equation}
\label{plustwee}
-
{2^{e+1}}\sum_{k=e+2}^{\infty}{1\over 2^k}\left({1\over [L(\zeta_{2^k}):\mathbb Q]}
-{1\over [L(\zeta_{2^{k+1}}):\mathbb Q]}\right),
\end{equation}
if $\epsilon=-1$.
\end{Lem}
{\it Proof}. This easily follows on combining the previous lemma 
with equation (\ref{heurieen}), respectively (\ref{heuritwee}). \qed\\

\noindent Remark. From (\ref{pluseen}) and (\ref{plustwee}) we infer that
$$\delta_2(-|r|)-\delta_2(|r|)=1-{3\over 2^{e+1}}+{2\over [L(\zeta_{2^{e+1}}):\mathbb Q]}
-{2\over [L(\zeta_{2^{e+2}}):\mathbb Q]}.$$

\noindent The number $\delta_2(r)$ can be readily evaluated on using the following simple fact
from algebraic number theory:
\begin{Lem}
\label{graad}
Let $K$ be a real quadratic field.
Let $k\ge 1$. Then $$[K(\zeta_{2^k}):\mathbb Q]=\cases{2^k &if $k\le 2$ or $K\ne \mathbb Q(\sqrt{2})$;\cr
2^{k-1} & if $k\ge 3$ and $K=\mathbb Q(\sqrt{2})$.}$$
\end{Lem}
{\it Proof}. If $K$ is a quadratic field other than $\mathbb Q(\sqrt{2})$ then there is an odd prime
that ramifies in it. This prime, however, does not ramify in $\mathbb Q(\zeta_{2^n})$, so in this
case $K$ and $\mathbb Q(\sqrt{2})$ are linearly disjoint.
Note that $\zeta_8+\zeta_8^{-1}=\sqrt{2}$ and hence $\mathbb Q(\sqrt{2})\subset \mathbb Q(\zeta_8)$.
Using the well-known result that $[\mathbb Q(\zeta_n):\mathbb Q]=\varphi(n)$, the result is then
easily completed. \qed\\

\noindent The result of this evaluation is stated below.
\begin{Lem}
We have $\delta_2(r)=\delta(r)$.
\end{Lem}

\noindent After all this preliminary work, it is straightforward to prove the two main results of this note:\\

\noindent {\it Proof of Theorem} \ref{maintwee}. 1) Left to the reader. 2) Combine the latter
lemma with Lemma \ref{bijna}. Comparison with Theorem \ref{oud} shows that
$H_{a,b}^{(2)}(x)\sim N_{a,b}(x)$ as $x\rightarrow \infty$ and thus $H_{a,b}^{(2)}(x)$ is
an asymptotically exact approximation of $N_{a,b}(x)$. \qed\\

\noindent {\it Proof of Theorem} \ref{main}. Combine part 2 of Theorem \ref{maintwee} 
(with any $A>3$), Theorem \ref{oud}
and equations (\ref{heurieen}) and (\ref{heuritwee}). \qed

\section{Two alternative formulations}
\subsection{An alternative formulation using Ramanujan sums}
Recall that the {\it Ramanujan sum} $c_n(m)$ is  defined as $\sum_{1\le k\le n,~(k,n)=1}e^{2\pi i km/n}$.
It is well-known that $c_n(m)\in \mathbb Z$ and, more in particular, that
$$c_n(m)=\varphi(n){\mu(n/(n,m))\over \varphi(n/(n,m))}.$$ This is known as 
{\it H\"older's identity}. It implies
that $c_n(m)=c_n((n,m))$. For our purposes the following weak version of H\"older's identity will
suffice:
\begin{equation}
\label{zwak}
c_{2^v}(t)=\cases{0 & if $\nu_2(t)<v$;\cr 
-\varphi(2^v) &if $\nu_2(t)=v-1$;\cr
\varphi(2^v) &if $\nu_2(t)\ge v$.}
\end{equation}
Another elementary property of Ramanujan sums we need is that for arbitrary natural numbers $n$ and $m$
\begin{equation}
\label{indicator}
{1\over n}\sum_{d|n}c_d(m)=\cases{1 &if $n|m$;\cr 0 & otherwise.}
\end{equation}
Suppose that $\nu_p(r)=0$, then ord$_r(p)[\mathbb F_p^*:\langle r \rangle]=p-1$. Note that ord$_r(p)$
is off iff $2^{\nu_2(p-1)}|[\mathbb F_p^*:\langle r\rangle]$. Using identity (\ref{indicator}) it then
follows that
\begin{equation}
\label{schietop}
N_{a,b}(x)=\pi(x)-\sum_{p\le x,~p\nmid 2ab}2^{-\nu_2(p-1)}\sum_{v\le \nu_2(p-1)}c_{2^v}
([\mathbb F_p^*:\langle r\rangle])+O(\omega(ab)).
\end{equation}
Corollary \ref{puhpuhpuh} below shows that if in the latter double sum the summation 
is restricted to those $v$ satisfying in
addition $v\le e$, respectively $v\le e+1$, then $K_{a,b}^{(1)}(x)$, respectively $K_{a,b}^{(2)}(x)$
is obtained. This in combination with Theorems \ref{oud}
and \ref{maintwee} leads to the following theorem:
\begin{Thm}
\label{maindrie}
We have, in the notation of Theorem {\rm \ref{main}},
$N_{a,b}(x)=$
$$\pi(x)-\sum_{p\le x,~p\nmid 2ab}2^{-\nu_2(p-1)}\sum_{2^v|(p-1,2h)}c_{2^{\nu}}
([\mathbb F_p^*:\langle {r}\rangle])
+O\left({x(\log \log x)^4\over \log^3 x}\right),$$
and 
$$\sum_{p\le x,~p\nmid 2ab}2^{-\nu_2(p-1)}\sum_{e+2\le v\le \nu_2(p-1)}c_{2^{\nu}}
([\mathbb F_p^*:\langle {r}\rangle])=O\left({x(\log \log x)^4\over \log^3 x}\right),$$
where the implied constant depends at most on $a$ and $b$.
\end{Thm}
Remark. Note that $v\le {\rm min}(\nu_2(p-1),e+1)$ is equivalent with $2^v|(p-1,2h)$.
\begin{Lem} 
\label{puhpuh}
Let $a,b,\epsilon$ and $e$ be as in Theorem {\rm \ref{main}} and 
let $p\nmid 2ab$.\\
{\rm 1)} We have
$$2^{-\nu_2(p-1)}\sum_{v\le {\rm min}(\nu_2(p-1),e)}c_{2^{\nu}}([\mathbb F_p^*:\langle {r}\rangle])
=k_{a,b}^{(1)}(p).$$
{\rm 2)} We have $$2^{-\nu_2(p-1)}\sum_{v\le {\rm min}(\nu_2(p-1),e+1)}c_{2^{\nu}}([\mathbb F_p^*:\langle 
{r}\rangle])
=k_{a,b}^{(2)}(p).$$
\end{Lem}
\begin{Cor} 
\label{puhpuhpuh}
For $1\le j\le 2$ we have
$$\sum_{p\le x,~p\nmid 2ab}2^{-\nu_2(p-1)}
\sum_{v\le {\rm min}(\nu_2(p-1),e+j-1)}c_{2^{\nu}}([\mathbb F_p^*:\langle {r}\rangle])=K_{a,b}^{(j)}(x).$$
\end{Cor}
{\it Proof of Lemma} \ref{puhpuh}. 1) Let us consider the case $\epsilon=-1$ and $\nu_2(p-1)>e$ (the remaining cases
are similar and left to the reader). Since $(-1)^{(p-1)/2^e}\equiv 1({\rm mod~}p)$ we see
that $-1$ and hence $r$ is a $2^e$th-power mod $p$ and thus 
$\nu_2([\mathbb F_p^*:\langle {r}\rangle])\ge e$. Hence the sum in the statement of the lemma reduces
to $2^{-\nu_2(p-1)}\sum_{v\le e}\varphi(2^v)=
2^{e-\nu_2(p-1)}=k_{a,b}^{(1)}(p),$
where (\ref{zwak}), (\ref{lokaaleen}) and the identity $\sum_{d|n}\varphi(d)=n$ are used.\\
2) The case $\nu_2(p-1)\le e$. The quantity under consideration agrees with that of part 1 and by 
(\ref{lokaaltwee}) we obtain that $k_{a,b}^{(1)}(p)=(1+\epsilon)/2=k_{a,b}^{(2)}(p)$.\\
The case $\nu_2(p-1)=e+1$. Now 
$$(-1)^{p-1\over 2^{e+1}}\equiv 1({\rm mod~}p),(r_0^h)^{p-1\over 2^{e+1}}\equiv 
({r_0\over p})({\rm mod~}p){\rm ~and~hence~}r^{p-1\over 2^{e+1}}\equiv \epsilon({r_0\over p})({\rm mod~}p).$$
It follows that
$\nu_2([\mathbb F_p^*:\langle r\rangle)\ge e+1$ if $\epsilon({r_0\over p})=1$ and
$\nu_2([\mathbb F_p^*:\langle r\rangle])=e$ if $\epsilon({r_0\over p})=-1$. Using (\ref{zwak}) the quantity
under consideration is seen to reduce to
$$2^{-\nu_2(p-1)}\left(\sum_{v\le e}\varphi(2^v)+\epsilon({r_0\over p})2^e\right)={1+\epsilon({r_0\over p})
\over 2}.$$
By (\ref{lokaaltwee}) this equals $k_{a,b}^{(2)}(p)$.\\
The case  $\nu_2(p-1)>e+1$. Now $r^{(p-1)/2^{1+e}}\equiv ({r_0\over p}){\rm mod~}p)$. Proceeding as before
the quantity under consideration reduces to
$$2^{-\nu_2(p-1)}\left(\sum_{v\le e}\varphi(2^v)+({r_0\over p})2^e\right)=2^{e-\nu_2(p-1)}(1+({r_0\over p})).$$

\subsection{An alternative formulation involving character sums}
Let $G$ be a cyclic group of order $n$ and $g\in G$. It is not difficult to show that, for any $d|n$,
$\sum_{{\rm ord}(\chi)=d}\chi(g)=c_d([G:\langle g\rangle])$. Using this and noting
that $\chi(r)=\chi(\epsilon)\chi^h(r_0)$, equation (\ref{schietop}) can be rewritten as
\begin{equation}
\label{schietop2}
N_{a,b}(x)=\pi(x)-\sum_{p\le x,~p\nmid 2ab}2^{-\nu_2(p-1)}\sum_{{\rm ord}(\chi)|2^{\nu_2(p-1)}}\chi(\epsilon)\chi^h(r_0)
+O(\omega(ab)),
\end{equation}
where the sum is over all characters of $\mathbb F_p^*$ having order dividing $2^{\nu_2(p-1)}$.
Note that if $\chi$ is of order $2^v$, then $\chi^h$ is the trivial character  if $v\le e$ and a
quadratic character if $v=e+1$. If in the main term of (\ref{schietop2}) only those characters of order dividing
$h$ are retained, i.e. those for which $\chi^h$ is 
the trivial character, then $H_{a,b}^{(1)}(x)$ is obtained (this is a reformulation of
part 1 of Lemma \ref{puhpuh}) and hence, by part 1 of 
Theorem \ref{maintwee} the na\"{\i}ve heuristic. If in (\ref{schietop2}) only those 
characters of order dividing
$2h$ are retained, i.e. those for which $\chi^h$ is 
the trivial or a quadratic character, then the 
asymptotically exact heuristic is obtained. 
The error term assertion in Theorem \ref{maindrie} can be reformulated as:
\begin{Prop}
We have
$$\sum_{p\le x,~p\nmid 2ab}2^{-\nu_2(p-1)}\sum_{2^{e+2}|{\rm ord}(\chi)|2^{\nu_2(p-1)}}
\chi(\epsilon)\chi^h(r_0)=O\left({x(\log \log x)^4\over \log^3 x}\right),$$
where the implied constant depends at most on $a$ and $b$.
\end{Prop}
In the setting of near primitive roots it is already known
that for the main term of the counting function of (near) primitive roots only the contributions 
coming from characters that are either trivial or quadratic need to be included \cite{M2}. 
\section{Conclusion}
There is a na\"{\i}ve heuristic for $N_{a,b}(x)$ that in many, but not all, cases is asymptotically exact.
There is a quadratic modification of this heuristic involving the Legendre symboll that is
{\it always} asymptotically exact. The same phenomenon is observed (assuming GRH) in the setting of Artin's 
primitive root conjecture.


\medskip\noindent {\footnotesize Korteweg-de Vries Institute, University of Amsterdam,
Plantage Muidergracht 24, 1018 TV Amsterdam, The Netherlands.\\
e-mail: moree@science.uva.nl }

\end{document}